\documentclass[12pt,a4paper]{article}
\usepackage{graphicx}
\usepackage{graphics}
\usepackage[usenames]{color}
\usepackage{colortbl}

\newcommand{\la}{\lambda}
\newcommand{\vp}{{\mathbf p}}
\newcommand{\vz}{{\mathbf z}}

\newcommand{\vy}{{\mathbf y}}
\newcommand{\vw}{{\mathbf w}}

\newcounter{theoremcounter}

\newtheorem{theorem}[theoremcounter]{Theorem}

\begin{document}

%%%%% TITLE %%%%%

\title{Bounds on the rate of convergence  for  Markovian queueing models with catastrophes} %in capital

\date{}
%%%%% AUTHORS %%%%%

\author{A. Zeifman\footnote{Vologda State University; Federal Research Center “Computer Science and Control” of RAS;
Vologda Research Center of RAS, Moscow Center for Fundamental and Applied Mathematics, Moscow State University, Russia; $a\_$zeifman@mail.ru}}

\maketitle

%\renewcommand{\baselinestretch}{1}

%\begin{document}

\medskip

{\bf Abstract.} In this note, a general approach to the study of non-stationary Markov chains with catastrophes and the corresponding queuing models is considered, as well as to obtain estimates of the limiting regime itself.
As an illustration, an example of a queuing model is studied.

\medskip

\medskip

\section{Introduction}

We consider a  general nonstationary Markovian queueing model under additional assumption of possibility of 
catastrophes of the system. As a rule this assumption is sufficient for ergodicity of the corresponding queue-length process.

There is a number of  investigations in this area, see, for instance, \cite{Ammar2014,Ammar2020,Chakravarthy2017,Chen1997,Chen2004,Di
Crescenzo2008,Dudin2001,Li2017,Zhang2015,Zeifman2012,Zeifman2017,Zeifman2017ecms,Zeifman2020} and references therein.  

 In these papers as a rule stationary distributions or transient behavior are studied. In our previous papers we  obtained estimates on the rate of convergence to the limiting regime for a number of classes of Markovian queueing models with catastrophes. This note is devoted to a simple and general method  for study of ergodicity of such models (in particular, in nonstationary situations). This approach enables us to more efficiently compute the main probabilistic characteristics for Markovian queueing models, as shown in \cite {Zeifman2021_mdpi}.

 Here we obtain upper bounds on the rate of convergence for such models and apply these estimates to some specific situations.

Let

$\la_{i,i+k}(t)$  be the intensity of arrival of group of $k$ customers to the queue at the moment $t$, if the current length of queue equals  $i$;

$\mu_{i,i-k}(t)$  be the intensity of service of a group of $k$ customers to the queue at the moment $t$, if the current length of  queue equals  $i$.

\smallskip

 In addition, we separately introduce a special notation  for the catastrophe (disaster) intensity, that is, the intensity of  simultaneous loss of all customers. Namely,  let 
$\beta_k(t)$  be a disaster (catastrophe) intensity, if the current size of the length of queue equals  $k$.

\medskip

Consider the corresponding  queue-length process $X(t)$. Then the
intensity matrix $Q(t)=\left(q_{ij}(t)\right)_{i,j=0}^{\infty}$ for
$X(t)$ takes the following form:

$$
Q\left(t\right)=\left(
\begin{array}{cccccccc}
q_{00}\left(t\right) & \la_{01}\left(t\right)  & \la_{02}\left(t\right)   & \la_{03}\left(t\right)  & \la_{04}\left(t\right) & \ldots  & \ldots\\[3pt]
\mu_{10}(t)+\beta_1\left(t\right)& q_{11}\left(t\right)  & \la_{12}\left(t\right)  & \la_{13}\left(t\right)   &  \ldots & \ldots  & \ldots\\[3pt]
\mu_{20}(t)+ \beta_2\left(t\right)  & \mu_{21}(t)\left(t\right)    & q_{22}\left(t\right)& \la_{23}(t)\left(t\right)  &   \la_{24}\left(t\right)  &  \ldots &   \ldots\\[3pt]
\ldots&\ldots&\ldots&\ldots&\ldots&\ldots&\ldots\\[3pt]
\mu_{j0}(t)+\beta_j\left(t\right) &  &  \ldots & \mu_{j,j-1}\left(t\right)  &   q_{jj}\left(t\right) & \la_{j,j+1}\left(t\right)  & \ldots\\[3pt]
\vdots&\vdots&\vdots&\vdots&\vdots&\vdots&\ddots
\end{array}
\right),
$$
\noindent where $q_{ii}\left(t\right)$'s are such that all row sums of
the matrix equal to zero for any $t \ge 0$.

\bigskip

 Applying the standard approach
(see for instance \cite{Zeifman2017,Zeifman2020porto}) we assume that
all the intensity functions $q_{ij}(t)$ are locally integrable on
$[0,\infty)$, and  that  $\sup_i|q_{ii}(t)| = L(t) < \infty,$ for almost all $t \ge 0$.

\bigskip

Then the probabilistic dynamics of the process $\{ X(t), \ t\geq
0\}$ is given by the forward Kolmogorov system
\begin{equation} \label{ur01}
\frac{d {\bf p}(t)}{dt}=A(t){\bf p}(t),
\end{equation}
\noindent where
$$
A(t) = Q^T(t)=\left(
\begin{array}{cccccccc}
q_{00}\left(t\right) & \mu_{10}\left(t\right)+\beta_1\left(t\right)    &  \mu_{20}\left(t\right)+\beta_2\left(t\right) &  \ldots  &  \mu_{j0}\left(t\right)+\beta_j\left(t\right) & \ldots\\[3pt]
\la_{01}\left(t\right)   & q_{11}\left(t\right)  &  \mu_{21}\left(t\right)   &  \mu_{31}\left(t\right)  &  \ldots & \ldots  & \ldots\\[3pt]
\la_{02}\left(t\right)  &   \la_{12}\left(t\right)  & q_{22}\left(t\right)&  \mu_{32}\left(t\right)  &    \mu_{42}\left(t\right)   &  \ldots &   \ldots\\[3pt]
\vdots & \ldots&\ldots&\ldots&\ldots&\ldots&\ldots\\[3pt]
\vdots&\vdots&\vdots&\vdots&\vdots&\vdots&\ddots
\end{array}
\right)
$$
\noindent is the transposed intensity matrix and $\vp(t)$ is the
column vector of state probabilities, $\vp(t) = \left(p_0(t),
p_1(t), \dots\right)^T$.

\bigskip
Then, applying the modified combined approach of
\cite{Zeifman2017ecms} and \cite{Zeifman2018} we can obtain bounds
on the rate of convergence of the queue-length process to its
limiting characteristics and compute them. We separately consider
the important special cases.

\section{Basic Notions}

Denote by $p_{ij}(s,t)=P\left\{ X(t)=j\left| X(s)=i\right.
\right\}$, $i,j \ge 0, \;0\leq s\leq t$ the transition probabilities
of $X(t)$ and by  $p_i(t)=P \left\{ X(t) =i \right\}$  the
probability that $X(t)$  is in state $i$ at time $t$. Let ${\bf
p}(t) = \left(p_0(t), p_1(t), \dots\right)^T$ be probability
distribution vector at instant $t$.

Throughout the paper by $\|\cdot\|$  we denote  the $l_1$-norm, i.
e.  $\|{\vp(t)}\|=\sum_{k} |p_k(t)|$, and $\|A(t)\| =
\sup_{j} \sum_{i} |a_{ij}(t)|$. Let
$\Omega$ be a set of all stochastic vectors, i. e. $l_1$ vectors with
non-negative coordinates and unit norm. Hence we have $\|A(t)\| =
2\sup_{k} |q_{kk}(t)| - 2 L(t) < \infty$ for almost all ${t
\ge 0}$. Hence the operator function $A(t)$ from $l_1$ into itself
is bounded for almost all $t \ge 0$ and locally integrable on
${[0,\infty)}$. Therefore we can consider (\ref{ur01}) as a
differential equation in the space $l_1$ with bounded operator.

It is well known (see \cite{Daleckij1974}) that the Cauchy problem
for differential equation (\ref{ur01}) has a unique solutions for an
arbitrary initial condition, and  $\vp(s) \in \Omega$ implies
$\vp(t) \in \Omega$ for ${t \ge s \ge 0}$.

Denote by $E(t,k) = E(X(t)|X(0)=k)$ the conditional expected number
of customers in the system at instant $t$, provided that initially
(at instant $t=0$) $k$ customers were present in the system.

Recall that a Markov chain \mbox{$\{ X(t), \ t\geq 0\}$} is called
{\it weakly ergodic}, if $\|\vp^{*}(t)-\vp^{**}(t)\| \to 0$ as $t
\to \infty$ for any initial conditions $\vp^{*}(0)$ and
$\vp^{**}(0)$, where $\vp^{*}(t)$ and $\vp^{**}(t)$ are the
corresponding solutions of (\ref{ur01});  and {\it exponentially  weakly ergodic} if the difference tends to zero exponentially fast. A Markov chain \mbox{$\{
X(t), \ t\geq 0\}$} has  the {limiting mean} $\varphi (t)$, if $
\lim_{t \to \infty }  \left(\varphi (t) - E(t,k)\right) = 0$ for any
$k$.

\section{Main Bounds}

Rewrite the forward Kolmogorov system (\ref{ur01}) as
\begin{equation}
\frac{d {\bf p} }{dt}=A^*\left( t\right) {{\bf p} }  +{\bf g}
\left(t\right), \quad t\ge 0. \label{eq112''}
\end{equation}

Here  ${\bf
g}\left(t\right)=\left(\beta_*\left(t\right),0,0, \dots\right)^T$,
$A^*\left(t\right)=\left (
a_{ij}^*\left(t\right)\right)_{i,j=0}^{\infty}$, and
\begin{eqnarray*}
a_{ij}^*\left(t\right) =   \left\{
\begin{array}{ccccccc}
q_{00}\left(t\right) - \beta_*\left(t\right), & \mbox { if }  i= j=0, \\
\mu_{j0}\left(t\right) +\beta_j(t)- \beta_*\left(t\right), & \mbox { if }  i= 0,\ j >0 \\
q_{ji}\left(t\right), & \mbox { otherwise },
\end{array}
\right.  \label{1101}
\end{eqnarray*}
\noindent where $\beta_*\left(t\right) = \inf_i
\beta_i\left(t\right)$.

\bigskip

Denote $\vy = \vp^{*}-\vp^{**}$.

Then
$$
\frac{d \vy(t)}{dt}= A^*(t)\vy(t). \label{2.0611}
$$

\smallskip

Let now $d_k$ be positive numbers for $k \ge 0$, and let
$$
d=\inf_k d_k >0, \quad d^*=\sup_k d_k  \le \infty. \label{2.0612}
$$

\smallskip

Put $w_k=d_ky_k$, for $k \ge 0$. Consider a new vector function $\vw(t)=D\vy(t)$, where 
$D$ is a  diagonal matrix with entries $d_k$.

\smallskip

Then we obtain
$$
\frac{d \vw(t)}{dt}= A^*_{{D}}(t)\vw(t), \label{2.06DD}
$$
\noindent where
$A^*_{{D}}(t)={D}A^*(t){{D}}^{-1}=\left(a_{i,j,{D}}^*(t)\right)_{i,j=0}^{\infty}$,
with the corresponding elements. 

\bigskip
Let 
\begin{equation}
\beta_{**}(t)= \inf_i\left(|a_{i,i,{D}}^*(t)| -\sum_{j\neq i} a_{j,i,{D}}^*(t)\right),
\label{lognorm2}
\end{equation}

  Then one can write the following estimate for the upper right derivative of $\|\vw(t)\|$
$$
\frac{d^+_r}{dt} \|\vw(t)\| \le -\beta_{**}(t) \|\vw(t)\|,
\label{lognorm1}
$$
\noindent and then, dividing by the $\|\vw(t)\|$ and integrating, one will have the following upper bound:
$$
\|\vw(t)\| \le e^{-\int_0^t \beta_{**}(u) \, du} \|\vw(0)\|. \label{lognorm3}
$$

\medskip

{\bf Remark.} In fact, there is also the usual right-hand derivative of the norm, this is the logarithmic norm, which we most often use (see for instance \cite{{Zeifman2021_mdpi}}), so its application would lead to the same result.

\bigskip

If we compare different norms of the vector, we get

$$
\| \vw(t)\|= \|D \vy(t)\|= \|D\left(\vp^{*}(t)-\vp^{**}(t)\right)\|, \label{norms01}
$$
\noindent and 
$$
d\| \vy(t)\| \le  \| \vw(t)\| \le d^*  \| \vy(t)\|.  \label{norms02}
$$

Hence we have the following statement.

\smallskip

\begin{theorem}  \hspace{-0.2cm}{\bf.}{} 
Let there exists a sequence $\{d_k, \  k \ge 0\}$ such that 
\begin{equation}
\label{bet*} \int_0^\infty \beta_{**}\left(t\right)\, dt = +\infty.
\end{equation}
Then the queue-length process $X\left(t\right)$ is weakly
ergodic and the following bound on the rate of convergence holds:
\begin{eqnarray}
\left\|{{\bf p}}^{*}\left(t\right)-{{\bf
p}}^{**}\left(t\right)\right\| \le d^{-1}
e^{-\int\limits_0^t \beta_{**}\left(\tau\right)\, d\tau}\left\|D\left({{\bf
p}}^{*}\left(0\right)-{{\bf p}}^{**}\left(0\right)\right)\right\|.
\label{20101}\end{eqnarray}
\smallskip
 Moreover,
 
\smallskip 

(i) if $d^*< \infty$ then $X\left(t\right)$ is weakly
ergodic in the uniform operator topology and the following bound
hold
\begin{eqnarray}
\left\|{{\bf p}}^{*}\left(t\right)-{{\bf
p}}^{**}\left(t\right)\right\| \le  \frac{2d^*}{d}
e^{-\int\limits_0^t \beta_{**}\left(\tau\right)\, d\tau},
\label{2011}\end{eqnarray} \noindent for any initial conditions
${{\bf p}}^{*}\left(0\right), {{\bf p}}^{**}\left(0\right)$ and any
$t \ge 0$.

\smallskip

(ii) if $d^*= \infty$, and in addition  $W=\inf_{i \ge 1} \frac {d_i}{i} > 0$, 
then  $X(t)$ has the limiting mean, say ${\rm\phi}(t)=E(t,0)$, and the
following bound holds:
\begin{equation}
|E(t,j) - E(t,0)| \le  \frac{d_0+d_j}{W}e^{-\int\limits_0^t
\beta_{**}(\tau)\, d\tau}, \label{3012}
\end{equation}
for any $j$ and any $t \ge 0$.
\end{theorem}

\bigskip

Let now there exist positive $R_{**}$ and $b_{**}$ such that 
\begin{equation}
e^{-\int\limits_s^t\beta_{**}(\tau)\, d\tau} \le R_{**}e^{-b_{**}(t-s)}, \label{3021}
\end{equation}
\noindent for any $0 \le s \le t$. Then $X(t)$ is exponentially weakly ergodic and we can estimate the limiting regime itself by the following way.

\smallskip

Let, in addition the 'common catastrophe rate' $\beta_*(t)$ be bounded, i.e. 
\begin{equation}
\beta_*(t)\le b^*< \infty    \mbox { for almost all } t \ge 0.\label{3022}
\end{equation}

\smallskip

Denote by $U(t,s)$ the Cauchy operator of equation (\ref{eq112''}), then the solution of this equation looks as
$$
\vp(t)=U(t,0)\vp(0)+ \int_0^t U(t,\tau){\bf g}(\tau)\, d \tau, \label{3023}
$$
\noindent where $\vp(0)$ is the initial probability distribution of $X(t)$ (initial condition).

Hence we have (in $1D$ norm, where $\|\vz\|_{1D} = \|D\vz\|$ and  $\|B\|_{1D}=\|DBD^{-1}\|_{1}$):
\begin{eqnarray*}
\|\vp(t)\|_{1D} \le \|U(t,0)\|_{1D} \|\vp(0)\|_{1D}+ \int_0^t \|U(t,\tau)\|_{1D}\|{\bf g}(\tau)\|_{1D} \, d \tau  \\ \le
R_{**}e^{-b_{**}t} \|\vp(0)\|_{1D} + \int_0^t R_{**}e^{-b_{**}(t-\tau)}d_0 b^* \, d \tau \le 
o(1)+ \frac{R_{**}d_0 b^*}{b_{**}}, \nonumber \label{3024}
\end{eqnarray*}
\noindent and the following statement.

\bigskip

\begin{theorem} \hspace{-0.2cm}{\bf.}{} 
Let for some positive sequence $\{d_k\}$ inequalities  (\ref{3021})  and (\ref{3022}) hold.

Then the existing by Theorem 1 limiting regime satisfied to the following bound:
\begin{equation}
\limsup_{t \to \infty}\|\vp(t)\|_{1D} \le  \frac{R_{**}d_0 b^*}{b_{**}}. \nonumber \label{3025}
\end{equation}

\end{theorem}

\bigskip

\bigskip

{\bf Example.} Consider here as an example the model from \cite{marin2020,Zeifman2020_marin} with additional disasters (catastrophes).

Then the corresponding intensity matrix $A(t)$ of $X(t)$ has the following form:

$$
A(t) = {\tiny \left(
\begin{array}{cccccccc}
-\lambda(t) & \mu(t)+\gamma_1(t) & \gamma_2(t) &  \gamma_3(t) & \dots \\
\la(t)b_1 & - \left (\lambda(t)B_2+ \mu(t) + \gamma_1(t) \right ) & \mu(t)   &  0 & \dots \\
\la(t)b_2& \la(t)b_2 & - \left ( \lambda(t)B_3 + \mu(t) + \gamma_2(t) \right ) & \mu(t)   &  \dots \\
\la(t)b_3 & \la(t)b_3 & \la(t)b_3  &  - \left ( \lambda(t)B_4 + \mu(t) + \gamma_3 (t)\right ) \dots \\
\vdots& \vdots& \vdots& \ddots \\
\end{array}
\right)}
$$
\noindent where $B_k=\sum_{i \ge k}b_i$, $B_1= 1$, all $b_k \ge 0$, and $\sum_k kB_k < \infty.$

\bigskip

In the previous studies the authors suppose geometric decreasing of $b_k$.  Here we outline the simple situation of 
slowly decreasing arrivals intensities, namely let $b_k=\frac{4}{k(k+1)(k+2)}$. Put $d_0=1$, $d_k=k+1$ for $k \ge 1$, and $\beta_*(t)=\inf_{i \ge 1} \gamma_i(t)$. Then one has in (\ref{lognorm2}) 
$$
\beta_{**}(t)= \beta_{*}(t)- \lambda(t)\sum_{k \ge 1}\left(d_k-1\right)b_k=\beta_{*}(t) - \frac{1}{2}\lambda(t).$$

In particular, if $\lambda(t)=2+2\cos 2 \pi t$, $\gamma_k(t)=2+\frac{1+\sin2 \pi t}{k}$, and $\mu(t)$ is an arbitrary $1-$periodic function,  that is, a function periodic in $t$ with period equals to 1. Then $\beta_{**}(t) =  1+2\cos 2 \pi t$.

Hence we have in $R_{**} \le 2$, $b_{**} = 1$, $b^* = 4$ in (\ref{3021}) and (\ref{3022}) respectively. Moreover $W=1$,
$X(t)$ is exponentially ergodic and has the limiting mean.  Theorems 1 and 2 give us the following bounds:

$$
\left\|{{\bf p}}^{*}\left(t\right)-{{\bf
p}}^{**}\left(t\right)\right\| \le 2 e^{-t}\left\|D\left({{\bf
p}}^{*}\left(0\right)-{{\bf p}}^{**}\left(0\right)\right)\right\|,
$$ 
$$
|E(t,j) - E(t,0)| \le  2\left(1+j\right)e^{-t}, 
$$
\noindent and
$$
\limsup_{t \to \infty}\|\vp(t)\|_{1D} = \limsup_{t \to \infty}E(t,0) \le  8.
$$

\bigskip

{\bf Acknowledgement.} This research was supported by Russian Science Foundation under grant 19-11-00020.
Author thanks the referees for helpful remarks.

\end{document}